\documentclass[11pt,twoside]{amsart}
\usepackage{amsmath, amssymb, amscd}
\def\proof@headerfont{\upshape\bfseries}
\theoremstyle{plain}
\newtheorem{thm}{Theorem}[section]
\newtheorem{cor}[thm]{Corollary}
\newtheorem{lem}[thm]{Lemma}
\newtheorem{prop}[thm]{Proposition}
\newcommand{\ds}{\displaystyle}
\newtheorem{defn}[thm]{Definition}
\newtheorem{rems}[thm]{Remarks}
\newtheorem{prel*}{Preliminaries}
\newtheorem{examples*}{Examples}

\newcommand{\D}{{\mathbb D}}

\newcommand{\Index}{\operatorname{index}}
 
\newcommand{\tr}{\operatorname{tr}}

\newcommand{\R}{\operatorname{\mathbb R}}
\newcommand{\Z}{\operatorname{\mathbb Z}}

\newcommand{\nc}{\newcommand}
\nc{\nt}{\newtheorem}
\nc{\gf}[2]{\genfrac{}{}{0pt}{}{#1}{#2}} \nc{\mb}[1]{{\mbox{$ #1 $}}}
\nc{\real}{{\mathbb R}}
\nc{\comp}{{\mathbb C}}
\nc{\ints}{{\mathbb Z}}
\nc{\Ltoo}{\mb{L^2({\mathbf H})}}
\nc{\rtoo}{\mb{{\mathbf R}^2}}
\nc{\slr}{{\mathbf {PSL}}(2,\real)}
\nc{\slz}{{\mathbf {SL}}(2,\ints)}
\nc{\su}{{\mathbf {PSU}}(1,1)}
\nc{\so}{{\mathbf {SO}}}
\nc{\hyp}{{\mathbb H}}
\nc{\disc}{{\mathbf D}}
\nc{\torus}{{\mathbb T}}

\nc{\ca}{{\mathcal A}}
\nc{\cag}{{{\mathcal A}^\Gamma}}
\nc{\cg}{{\mathcal G}}
\nc{\chh}{{\mathcal H}}
\nc{\ck}{{\mathcal B}}
\nc{\cl}{{\mathcal L}}
\nc{\cm}{{\mathcal M}}
\nc{\cs}{{\mathcal S}}
\nc{\cz}{{\mathcal Z}}
\nc{\sind}{\sigma{\rm -ind}}

\begin{document}

\title{Quantum Hall Effect and Noncommutative Geometry}

\author{A.L. Carey, K.C. Hannabuss, V. Mathai}
\date{}
\maketitle

%\comm{Communicated by XXX}

\begin{abstract} We study magnetic Schr\"odinger operators with
random or almost periodic electric potentials on the hyperbolic plane,
motivated by the quantum Hall effect (QHE)
in which the hyperbolic geometry provides an effective
Hamiltonian. In addition we add some refinements to earlier results. 
We derive an analogue of the Connes-Kubo formula for the Hall
conductance via the quantum adiabatic theorem, 
identifying it as a geometric invariant associated to an algebra
of observables that turns out to be a crossed product algebra.
 We modify the Fredholm modules defined in \cite{CHMM} in order to prove
the integrality of the Hall conductance in this case.
\end{abstract}

%%%%%%%%%%%%%%%%%

\section*{Introduction}

In \cite{CHMM}, continuous and discrete
magnetic Hamiltonians containing terms arising from
a background hyperbolic geometry
were introduced. These may be thought of as effective
Hamiltonians for an analogue of the quantum Hall
effect studied in a
Euclidean model by Bellissard \cite{Bel+E+S} and Xia \cite{Xia}.
We interpret these Hamiltonians, following a suggestion of
Bellissard, as modelling spinless  electrons in a conducting material
with a perturbation term arising from a background
hyperbolic geometry. (In \cite{CHMM} we took the somewhat different view
that the conducting material exhibited hyperbolic geometry.)
They motivate constructing Fredholm modules
associated in a natural way with Riemann surfaces and
two dimensional orbifolds which give a higher genus
analogue of the work of Bellissard (which is the genus one case)
on the quantum Hall effect. In \cite{CHMM} we considered
Hamiltonians invariant under a projective action of a
Fuchsian group $\Gamma$.
We will only discuss groups whose actions
on hyperbolic space are free here
and refer the reader to \cite{MM} for the more general case.
In this paper we allow in addition a random potential
(which may be thought of as modelling impurities)
so that the invariance of the Hamiltonian under the Fuchsian group is replaced
by a type of ergodicity assumption.
There is an analogue of the Connes-Kubo cocycle of the Euclidean case for
the Hall conductance. This cocycle takes values which
are integer multiples of a fundamental unit in the case of free
actions and rational multiples for non-free actions.
Integrality follows by showing that
the cocycle gives the index of a certain Fredholm operator
(the conductance may also be thought of in terms of a topological index).
    Thus the models
in \cite{CHMM, MM}
fit the  noncommutative geometry
framework for magnetic Hamiltonians (see \cite {Co2}).

We begin by reviewing the construction of
magnetic Hamiltonians in a continuous
model with a background hyperbolic geometry term. There are also
discrete versions which are generalised Harper operators
\cite{Sun,CHMM,CHM,MM}. Our model of
   hyperbolic space is the upper half-plane $\hyp$ in
$\comp$ equipped with its usual Poincar\'e metric $(dx^2+dy^2)/y^2$, and
symplectic area form $\omega_\hyp = dx\wedge dy/y^2$. The group $\slr$ acts
transitively on $\hyp$ by M\"obius transformations $$ x+iy = \zeta \mapsto
g\zeta = \frac{a\zeta+b}{c\zeta+d},\quad\mbox{for } g=\left(
\begin{array}{cc} a & b\\ c & d
\end{array}\right).$$ Any Riemann surface of genus $g$ greater than 1 can
be realised as the quotient of
$\hyp$ by the action of its fundamental group realised as a cocompact 
torsion-free discrete subgroup $\Gamma$ of $ \slr$.

Pick a 1-form $\eta$ such that $d\eta = \theta\omega_\hyp$, for some fixed
$\theta \in \R$. As in geometric quantisation we may regard $\eta$ as
defining a connection $\nabla = d-i\eta$ on a line bundle $\cl$ over
$\hyp$, whose curvature is $\theta\omega_\hyp$. Physically we can think of
$\eta$ as the electromagnetic vector potential for a uniform magnetic field
of strength $\theta$ normal to $\hyp$. Using the Riemannian metric the
Hamiltonian of an electron in this field is given in suitable units by $$H
= H_\eta = \frac 12\nabla^*\nabla = \frac 12(d-i\eta)^*(d-i\eta).$$ Comtet
\cite{Com} has shown that $H$ differs from a multiple of the Casimir element
for $\slr$, $\frac 18{\bf J}.{\bf J}$, %+\frac 14B^2$ by a constant, where
$J_1$, $J_2$ and $J_3$ denote a certain representation of generators of the
Lie algebra $sl(2,\real)$, satisfying $[J_1,J_2] = -iJ_3,$  $[J_2,J_3]
= iJ_1,$ $[J_3,J_1] = iJ_2,$ so that ${\bf J}.{\bf J} =
J_1^2+J_2^2-J_3^2$ is the quadratic Casimir element showing
the underlying $\slr$-invariance of the theory. Comtet has computed
the spectrum of the unperturbed Hamiltonian $H_\eta$, for $\eta = -\theta
dx/y$, to be the union of finitely many eigenvalues $\{(2k+1)\theta
-k(k+1):k=0,1,2\ldots < \theta-\frac 12\}$, and the continuous spectrum
$[\frac 14 + \theta^2, \infty)$. Any $\eta$ is cohomologous to $-\theta
dx/y$ (since they both have $\omega_\hyp$ as differential) and forms
differing by an exact form $d\phi$ give equivalent models: in fact,
multiplying the wave functions by $\exp(i\phi)$ shows that the models for
$\eta$ and $-\theta dx/y$ are unitarily equivalent. This equivalence also
intertwines the $\Gamma$-actions so that the spectral densities for the two
models also coincide.

This Hamltonian can
be perturbed by adding a potential term $V$.
In \cite{CHMM}, we took
$V$ to be invariant under $\Gamma$. In \cite{CHM}
we allowed any smooth random potential function $V$ on $\mathbb H$
using two general notions of random potential
(in the literature random usually refers to
the $\Gamma$-action on the disorder space
being required to admit an ergodic
invariant measure). The class of random potentials we consider
here contains
    any smooth bounded potential $V$.
The perturbed Hamiltonian $H_{\eta,V} = H_\eta +V$
has unknown spectrum for such general $V$. However we are
able to deduce some qualitative aspects of the spectrum of these
Hamitonians by using a reduction (via Morita equivalence) to a simpler
case: that of a discrete model.

In Section 2, we extend the hyperbolic Connes-Kubo formula for the Hall
conductance for the continuous model in \cite{CHMM},
to the non-periodic case. We show that this hyperbolic Connes-Kubo
cocycle is cohomologous to another cyclic 2-cocycle
which is the Chern character of a Fredholm module, 
from which we can deduce that
the Hall conductance takes on integral values in
$2(g-1)\Z$ ($g>1$ being the genus). This result has been generalized in
\cite{MM} where for general cocompact Fuchsian groups $\Gamma$, it is shown
that the conductance takes on values in $\phi\Z$, where $\phi$ denotes
the orbifold Euler characteristic of the orbifold $\hyp/\Gamma$, i.e.
the conductance can take on certain fractional values.
In the Appendix we give a derivation of the hyperbolic Connes-Kubo 
formula for the Hall conductance, using the quantum adiabatic theorem
and standard physical reasoning.

%%%%%%%%%%%%%%%%%%%%%%%%
\section{Continuous model}
%%%%%%%%%%%%%%%%%%%%%%%%

\subsection{The geometry of the hyperbolic plane}
The upper half-plane can be mapped by the Cayley transform $z =
(\zeta-i)/(\zeta+i)$ to the unit disc $\D$ equipped with the metric
$|dz|^2/(1-|z|^2)^2$ and symplectic form $dz\,d\overline{z}/2i(1-|z|^2)^2$,
on which $\su$ acts, and some calculations are more easily done in that
setting. In order to preserve flexibility we shall work more abstractly
with a Lie group $G$ acting transitively on a space $X \sim G/K$. Although
we shall ultimately be interested in the case of $G = \slr$ or $\su$, and
$K$ the maximal compact subgroup which stabilises $\zeta = i$ or $z=0$ so
that $X=\mathbb H$ or $X= \mathbb D$, those details will play little role
in many of our calculations, though we shall need to assume that $X$ has a
$G$-invariant Riemannian metric and symplectic form $\omega_X$. We shall
denote by $\Gamma$ a discrete subgroup of $G$ which acts freely on $X$ and
hence intersects $K$ trivially.

We shall assume that $\cl$ is a hermitian line bundle over $X$, with a
connection, $\nabla$, or equivalently, for each pair of points $w$ and $z$
in $X$, we denote by $\tau(z,w)$ the parallel transport operator along the
geodesic from $\cl_w$ to $\cl_z$. In $\hyp$ with the line bundle
trivialised and $\eta = \theta dx/y$ one can calculate explicitly that
$$\tau(z,w) = \exp\left(i\int_w^z\eta\right ) =
[(z-\overline{w})/(w-\overline{z})]^\theta.$$ For general $\eta$ we have
$\eta - \theta dx/y = d\phi$ and $$\tau(z,w) = \exp(i\int_w^z\eta) =
[(z-\overline{w})/(w-\overline{z})]^\theta\exp(i(\phi(z)-\phi(w))).$$
Parallel transport round a geodesic triangle with vertices $z$, $w$, $v$,
gives rise to a holonomy factor:
$$\varpi(v,w,z) = \tau(v,z)^{-1}\tau(v,w)\tau(w,z),$$ and this is clearly
the same for any other choice of $\eta$, so we may as well work in the
general case.

\begin{lem} The holonomy can be written as
$\varpi(v,w,z) = \exp\left(i\theta\int_\Delta \omega_\hyp\right)$, where
$\Delta$ denotes the geodesic triangle with vertices $z$, $w$ and $v$. The
holonomy is invariant under the action of $G$, that is $\varpi(v,w,z) =
\varpi(gv,gw,gz)$, and under cyclic permutations of its arguments.
Transposition of any two vertices inverts $\varpi$. For any four points $u$
,$v$, $w$, $z$ in $X$ one has $$\varpi(u,v,w)\varpi(u,w,z) =
\varpi(u,v,z)\varpi(v,w,z).$$ 
\end{lem}

\subsection{Algebra of observables and random or almost periodic
potentials} The algebra of physical observables that we consider in the
continuous model should include the operators $f(H_{\eta, V})$ for any
bounded continuous function $f$ on $\R$ and for any smooth random potential
function
$V$ on $\mathbb H$ with disorder space $\Omega$.
We will see that the twisted $C^*$-algebra of the groupoid
$\mathcal G = \Gamma\backslash (X \times X \times \Omega)$, twisted by
$\varpi$, is large enough to contain all such operators. This
algebra also turns out to be the twisted $C^*$-algebra of the foliation
$\Omega_\Gamma$. This $C^*$-algebra is strongly
Morita equivalent to the cross product $C^*$-algebra $C(\Omega)
\rtimes_\sigma \Gamma$, where $\sigma$ is a multiplier on $\Gamma$ which is
determined by $\varpi$.

{\bf Assumptions} The {\em disorder space} $\Omega$ we assume to be
compact, to admit a Borel probability measure
$\Lambda$;
and that there is a continuous action of $\Gamma$ on $\Omega$ with a
dense orbit. 

The geometrical data described in the last subsection enables us to easily
describe the first of the two $C^*$algebras which appear in the theory.
This twisted algebra of kernels, which was introduced by Connes \cite{Co2}
is the $C^*$-algebra $\ck$ generated by compactly supported smooth
functions on
$X\times X \times \Omega$ with the multiplication $$k_1*k_2(z,w, r) =
\int_X k_1(z,v, r)k_2(v,w, r)\varpi(z,w,v)\,dv,$$ (where $dv$ is the
$G$-invariant measure defined by the metric) and $k^*(z,w, r)=
\overline{k(w,z, r)}$. The trace on $\ck$ is given by,
$\tau_\ck (k) = \int_{X\times \Omega} k(z,z, r)\,dz d\Lambda(r)$. Observe
that $X\times X\times\Omega$ is a groupoid with space of units $X\times
\Omega$ and with source and range maps $s((z,w,r)) = (w,r)$ and
$r((z,w,r')) = (z,r')$. Then the algebra of twisted kernels is the
extension of the $C^*$-algebra of the groupoid $X\times X\times\Omega$
defined by the cocycle $((v,w, r),(w,z, r)) \mapsto \varpi(v,w,z)$,
\cite{Ren1}.

\begin{lem} The algebra $\ck$ has a representation $\pi$ on the space
$\chh$ of
$L^2$ sections of $\cl \to X\times \Omega$ defined by
$$(\pi(k)\psi)(z, r) = \int_X k(z,w, r)\tau(z,w)\psi(w, r)\,dw.$$
\end{lem}

We now pick out a
$\Gamma$-invariant subalgebra $\ck^\Gamma$ of $\ck$. This condition reduces
simply to the requirement that the kernel satisfies
$k(\gamma^{-1}z,\gamma^{-1}w, \gamma^{-1}r) = k(z,w, r)$ for all $\gamma\in
\Gamma$. As before, observe that $\Gamma\backslash(X\times X\times\Omega)$
is a groupoid whose elements are $\Gamma$ orbits $(x,y, v)_\Gamma =
\{(\gamma x,\gamma y, \gamma v): \gamma \in \Gamma\}$ , with source and
range maps
$s((x,y, v)_\Gamma) = ( y, v)$ and $r((x,y, v)_\Gamma) = ( x, v)$. The
space of units is $\Omega_\Gamma = \Gamma\backslash(X\times \Omega)$. Then
the algebra of invariant twisted kernels $\ck^\Gamma$ is the extension of
the $C^*$-algebra of the groupoid $\Gamma\backslash(X\times X\times\Omega)$
defined by the cocycle $((v,w, r),(w,z, r)) \mapsto \varpi(v,w,z)$,
\cite{Ren1}. With our assumptions on
    the disorder space $\Omega$, there is in general {\em no} trace on the
algebra $\ck^\Gamma$, and there may not even be a weight on this algebra in
general. However, we mention that under the additional assumption that the
measure $\Lambda$ on $\Omega$ is $\Gamma$-invariant, the natural trace
$\tau_{\ck^\Gamma}$ for this algebra is given by the same formula as before
except that the integration is
    now over $\Omega_\Gamma = \Gamma \backslash (X\times \Omega)$ rather than
$X\times \Omega$, where we have identified $\Omega_\Gamma$ with a
fundamental domain: $ \tau_{\ck^\Gamma} (T)=\int_{\Omega_\Gamma} T(z,z,
r)dzd\Lambda(r). $ We also mention that under the additional assumption
that the measure $\Lambda$ on $\Omega$ is quasi-$\Gamma$-invariant, the
natural tracial weight $\tau_{\ck^\Gamma}$ for this algebra is given by $
\tau_{\ck^\Gamma} (T)=\int_{X\times\Omega} f(z,r)^2T(z,z, r)dzd\Lambda(r),
$ where $f\in C_c(X\times\Omega)$ is such that $\sum_{\gamma\in \Gamma}
(\gamma^*f)^2 = 1$.

We now recall a notion due to Connes \cite{Co2}.
\begin{defn} A {\em random} or {\em almost periodic} potential on $X$ is a
continuous family of smooth functions on the disorder space, $ \Omega \ni r
\mapsto V_r \in C^\infty(X)$ where the following equivariance is imposed:
$$ V_{\gamma r} = \gamma^* V_r \qquad \forall \gamma \in \Gamma, \forall
r\in \Omega.
$$
\end{defn}
\begin{rems} If $V$ is a $\Gamma$-invariant potential on $X$, then it is
clearly random for any disorder space. More generally, if $V$ is a
arbitrary smooth function on $X$ such that the set $\left\{ \gamma^* V :
\gamma \in \Gamma\right\}$ has compact closure in the strong operator
topology in $B(L^2(X))$, then $V$ is a random potential. 
\end{rems}
The reason the Hamiltonian can be accommodated within
the algebra $\ck^\Gamma$ is not hard to explain. Fix a base
point $u\in\D$ and introduce: $$\sigma(x,y) = \varpi(u,xu,xyu)$$
$$\phi(z,\gamma) = \varpi(u, \gamma^{-1}u,\gamma^{-1}z)\tau(u,z)^{-1}
\tau(u,\gamma^{-1}z).$$ Then $\sigma$ is the group 2-cocycle in the
projective action of $\su$ on $L^2(\D)$ defined by:
$$U(\gamma)\psi(z) = \phi(z,\gamma)\psi(\gamma^{-1}z)$$ where
$\psi\in L^2(\D), \gamma\in\su$. Note that $U$ is constructed so that the
$\Gamma$-invariant algebra
$\pi(\ck^\Gamma)$ is the intersection of $\pi(\ck)$ with the commutant of
$U$. Recall that the unperturbed Hamiltonian $H = H_\eta$ commutes with the
projective representation $U$ (cf. Lemma 4.9, \cite{CHMM}). So we see that
$H$ is affiliated to the von Neumann algebra generated by the
representation $\pi$ of $\ck^\Gamma$ (cf. Corollary 4.2 \cite{CHMM}).

A random potential $V$ can be viewed as defining an equivariant family of
Hamiltonians $\Omega \ni r \mapsto H_{\eta, V_r} = H + V_r \in {\rm
Oper}(L^2(X)$ where ${\rm Oper}(L^2(X))$ denotes closed operators on
$L^2(X)$.
Br\"uning and Sunada have proved an estimate on the Schwartz kernel of the
heat operator for any elliptic operator, and in particular for $\exp(-t
H_{\eta, V_r})$ for $t > 0$, which implies that it is $L^1$ in each
variable separately. Since this kernel is $\Gamma$-equivariant it follows
(in exactly the same fashion as Lemma 4 of \cite{BrSu}) that this estimate
implies that $\exp(-t H_{\eta, V_r})$ is actually in the algebra
$\ck^\Gamma$.

\begin{lem} One has
$f(H_{\eta, V}) \in \ck^\Gamma$ for any bounded continuous function $f$ on
$\R$ and for any random potential $V$ on $X$. In particular, the spectral
projections of $H_{\eta, V}$ corresponding to gaps in the spectrum lie in
$\ck^\Gamma$. 
\end{lem}

Following \cite{Bel+E+S},\cite{Nak+Bel} but using our weaker assumptions
we now have the
\begin{thm} Let $V$ be a smooth bounded function on $X$. Then $V$ is a
random potential for some disorder space $\Omega$ and therefore $f(H_{\eta,
V}) \in \ck^\Gamma$ for any bounded continuous function $f$ on $\R$.
\end{thm}
\noindent{\em Example}. Let the Iwasawa decomposition of $\su$ be written
$KAN$ then
$PSL(2,\Z)$ acts on $\D=\su/K$ by M\"obius transformations so that
$\Gamma\subset PSL(2,\Z)$ also acts. Let
$g_{\lambda,w}(z)=\lambda\frac{1-|z|^2}{|w-z|^2}$ where $\lambda\in
\R^+\cong A$, $w\in U(1)\cong K$ and $z\in\D$. Now let $\gamma=\left(
\begin{array}{cc}
\alpha & \beta\\
\bar\beta & \bar\alpha
\end{array}\right)$ and we calculate
$$ U(\gamma)g_{\lambda,w}U(\gamma^{-1})=g_{\lambda_{\gamma,w}\lambda,\gamma
w} $$ where $\lambda_{\gamma,w}= |\bar\beta w +\bar\alpha|^{-2}$. The
stabiliser of $g_{1,1}$ is
$\{\pm\left(
\begin{array}{cc} 1-in & in\\ -in & 1+in
\end{array}\right):\ n\in\R\}.$ This group is $MN$ where $MAN$ is the
maximal parabolic subgroup. Thus we have the usual action of $\su$ on
$\su/MN$ and hence {\em a fortiori} a $\Gamma$-action which is known to be
ergodic, cf. \cite{Zim}. Note that, regarding $\{e^{-g_{\lambda,w}}\}$ as a
set of bounded multiplication operators on $L^2(\D)$, the strong closure of
$\{U(\gamma)e^{-g_{\lambda,w}}U(\gamma^{-1})\ |\ \lambda\in\R, w\in U(1)\}$
is homeomorphic to $S^2$. (This is because taking the strong closure adds
the zero and identity operator to the set.) Thus in this example the
disorder space is $S^2$ which admits a dense orbit and a quasi-invariant
ergodic probability measure.

\subsection{Morita equivalence}
Our ability to calculate the possible
values of our generalised Connes-Kubo cocycle
rests on a Morita equivalence argument
due initially to \cite{M+R+W}. We use
the twisted version, \cite{Ren2}, \cite{Ren3}. We have already
noted that $\ck$ is the $C^*$-algebra of an extension of the groupoid
$X\times X\times \Omega$ by a cocycle defined by $\varpi$, and $\Gamma$
invariance of $\varpi$ means that $\ck^\Gamma$ is likewise the
$C^*$-algebra of an extension of $\Gamma\backslash(X\times X\times \Omega)$
by $\varpi$, where $\Gamma\backslash(X\times X\times \Omega)$ denotes the
groupoid obtained by factoring out the diagonal action of $\Gamma$. More
precisely, the groupoid elements are $\Gamma$ orbits $(x,y, v)_\Gamma =
\{(\gamma x,\gamma y, \gamma v): \gamma \in \Gamma\}$, with source and
range maps
$s((x,y, v)_\Gamma) = ( y, v)$ and $r((x,y, v)_\Gamma) = ( x, v)$.
Therefore $(x_1,y_1, v_1)_\Gamma$ and $(x_2,y_2, v_2)_\Gamma$ are
composable if and only if $y_1 = \gamma x_2$ and $v_1 = \gamma v_2$ for
some $\gamma \in \Gamma$, and then the composition is $(x_1,\gamma y_2,
\gamma v_2)_\Gamma$. We also note that $\Omega \times \Gamma$ is a
groupoid. The source and range maps are $s((v, \gamma)) = \gamma v$ and
$r((v, \gamma)) = v$. Therefore the elements
$(v_1, \gamma_1)$ and $(v_2, \gamma_2)$ are composable if and only if $v_1
= \gamma_2 v_2$, and the composition is $(\gamma_2^{-1} v_1, \gamma_1
\gamma_2)$.
\begin{thm} The algebra $\ck^\Gamma$ is Morita equivalent to the twisted
cross product algebra
$C(\Omega) \rtimes_{\bar\sigma} \Gamma$. \end{thm}
The proof is a consequence of:
\begin{lem} The line bundle $\cl$ over $X\times \Omega$ provides an
equivalence (in the sense of \cite{Ren2} Definition 5.3) between the
groupoid extensions $(\Gamma\backslash(X\times X\times \Omega))^\varpi$ of
$\Gamma\backslash(X\times X\times \Omega)$ defined by $\varpi$ and
$(\Omega \times\Gamma)^\sigma$ of $\Omega \times\Gamma$ defined by
$\overline{\sigma}$.
\end{lem}
Using the orientation reversing diffeomorphism of the
Riemann surface $\Sigma = \Gamma\backslash X$, one can show as in
Proposition 7 \cite{CHMM} that the algebra $C(\Omega)
\rtimes_{\overline{\sigma}} \Gamma$ is isomorphic to
$C(\Omega) \rtimes_{\sigma}\Gamma$, where $\bar\sigma$ denotes the complex
conjugate of $\sigma$.
Morita equivalence of algebras implies their $K$-groups are the same.
It is possible to
calculate the values taken by our cyclic cocycles
for the continous model in terms those taken by explicit cocycles
on $C(\Omega) \rtimes_{\bar\sigma} \Gamma$. The method
uses generalisations of arguments first developed for the
study of the Baum-Connes conjecture. Full details are in \cite {Ma}
and \cite {CHMM}.

\subsection{A hyperbolic Connes-Kubo formula, part I}
The quotient $\Sigma= \hyp/\Gamma$  is a Riemann
surface when $\Gamma$ is a cocompact torsion free subgroup of ${\bf PSL}(2, \mathbb R)$.
On a Riemann surface it is natural to investigate changes in the potential
corresponding to adding multiples of the real and imaginary parts of
holomorphic 1-forms.
(For the genus one case with an imaginary period this amounts to choosing
forms whose integral round one sort of cycle vanishes but the
    integral around the
other cycle is non-trivial.
Physically this
would correspond to putting a non-trivial voltage across one cycle and
measuring a current round the other.)

We let $a_j,j=1,2,\ldots,2g$ be a normalized symplectic basis of harmonic 1-forms 
on $\Sigma = \hyp/\Gamma$ where 
 $a_{j+g} = * a_j,\, j=1,2,\ldots,g,$
 and 
 $
 \int_\Sigma a_j \wedge a_{j+g} = 1
$ for all $j=1, \ldots, g.$ 
 % are dual to $B_j,j=1,2,\ldots,g$).
We introduce the map from $\hyp$ to
    $\real^{2g}$ given by $\Xi:z\mapsto 
(\int_u^z a_1,\ldots,\int_u^z a_{2g})$.
It is the lift to $\hyp$ of the Abel-Jacobi map,
    \cite{GH}
(this map is usually regarded as mapping from $\Sigma_g$ to the Jacobi variety
however we are thinking of it as a map between the universal
covers of these spaces).
Notice that $\Xi$ 
gives the period lattice in $\real^{2g}$ (that is the lattice determined
by the periods of the harmonic forms $a_j$) to be  the standard
integer lattice $\ints^{2g}$ so that $J(\Sigma_g)=\real^{2g}/\ints^{2g}$.
We give  $\real^{2g}$
the distinguished basis consisting of
the vertices in this integer period lattice.
We write for
the  corresponding coordinates $u_1,u_2,\ldots u_{2g}$. 
Let $\omega_J = \sum_{j=1}^g du_j \wedge du_{j+g}$ denote
the symplectic form on $\real^{2g}$. 
The closed 1-forms $c_j=\Xi^*(du_j)$ 
are cohomologous to $a_j$ for all $j=1, \ldots, 2g,$
and therefore we have

\begin{lem}\label{lem:Jac}
In the notation above,
$\Xi^*(\omega_J)$ is cohomologous 
to $ \sum_{j=1}^g a_j \wedge a_{j+g}.$
\end{lem}
%\begin{defn}
%Let
%$\delta_j\tau(z,w) =\displaystyle i\int_w^z c_j\,\tau(z,w)$
%and
%$$\delta_j\varpi(u,g^{-1}u,g^{-1}z) =
%i\int_{\partial\Delta} c_j\,\varpi(u,g^{-1}u,g^{-1}z)$$
%where $\Delta$ is a triangle with vertices at the three arguments of $\varpi$.
%\end{defn}

Suppose that $\alpha \in\ck$ is a kernel
    {\it decaying rapidly}. By this we mean that it
    satisfies an estimate
$$
|\alpha (x,y, r)| \le \phi(d(x,y)), \qquad r\in\Omega,
$$
where $\phi$ is a positive and rapidly decreasing function on $\mathbb R$.
Now define 
$$\delta_j\alpha = [\Omega_j,\alpha], \qquad \text{i.e.} \quad \delta_j\alpha (x,y, r) = 
(\Omega_j(x) - \Omega_j(y)) \alpha(x,y, r),$$
where $\Omega_j(z) =\displaystyle i\int_u^z a_j$. 
Since $\Omega_j(\gamma. z)  - \Omega_j(z)$ is a constant depending only on $\gamma$
but independent of $z$, and $|\Omega_j(\gamma. z)  - \Omega_j(z)| \le 
C ||a_j||_\infty d(z, \gamma.z) \le C_j \ell(\gamma)$, where $ ||a_j||_\infty$ is the 
supremum norm of $a_j$, 
$\gamma\in \Gamma$, $d(z, \gamma.z)$ is the Riemannian distance 
between $z$ and $\gamma z$, and $\ell(\gamma)$ is the word length of $\gamma$.
It follows that $\delta_j\alpha$ lies in $\ck$ and therefore
$\delta_j$ is a densely defined derivation
on the algebra $\ck$,  and hence also on  $ \ck^\Gamma$ since clearly if 
$\alpha$ is $\Gamma$-invariant, then so is $\delta_j \alpha$.

We may summarise the previous discussion
as

\begin{lem} For operators $A_0, A_1, A_2$ in
    $\ck^\Gamma$ whose integral kernels
are rapidly decaying we have cyclic cocycles defined by
$$c_{j,k}(A_0, A_1, A_2)=
\tr_{\ck^\Gamma}(A_0[\delta_jA_1,\delta_kA_2])
= \tr_{\ck^\Gamma}(A_0[\Omega_j, A_1][\Omega_k,A_2])$$
for $j,k=1,\ldots,2g$.
\end{lem}

The cyclic cocycle $c_{jk}$ can be interpreted as the Kubo
formula for the conductance due to
   currents in the $k$ direction induced by electric
fields in the $j$ direction, as explained in the Appendix.

%%%%%%%%%%%%%%%%%
\section{A Fredholm module}
%%%%%%%%%%%%%%%%%%

We shall now assume that $X$ has a spin structure, and we write $\cs$ for the
spin bundle.
The representation of $\ck^\Gamma$ on $\chh$ can then be extended to an
action on $\chh\otimes\cs$.
This module can be equipped with a Fredholm structure by taking $F$ to be
Clifford multiplication by a suitable unit vector (to be explained below), and
using the product of the trace on $\chh$ and the graded trace on
the Clifford algebra.
(If $\varepsilon$ denotes the grading operator on the spinors then the graded
trace is just $\tr\circ\varepsilon$.)

The same module can also be described more explicitly: it splits into
$\chh\otimes\cs^+\oplus \chh\otimes\cs^-$ (with the
    superscripted sign indicating the eigenvalue of $\varepsilon$).
Suppose that $\varphi$ is a $U(1)$ valued function on the group, which
satisfies
$\varphi(kgh) = \chi_1(k)\varphi(g)\chi_2(h)$ for $k$ and $h$ in $K$ and some
$\sigma$-characters $\chi_1$ and $\chi_2$ of $K$.
The involution $F$ can be taken to be the
    matrix multiplication operator:
$F = \left(\begin{array}{cc}
0 &\varphi^*\cr \varphi &0\cr
\end{array}\right).$
We may take for $\varphi$
the function used by Connes \cite{Co2}
which is essentially the Mishchenko element.
In the next subsection we will see that the module is 2-summable
for suitably decaying kernels.
Since
$\varphi$ is invariant under similtaneous conjugation of both variables
by elements of $\Gamma$,
$F$ preserves the $\Gamma$-invariant subspace.
%Let $X_\Gamma=\Gamma\backslash\hyp/K$.
%We then have:

\begin{thm}\label{thm:Fred} 
There is a dense subalgebra ${\mathcal B}^\Gamma_0$
of  ${\mathcal B}^\Gamma$ stable under the holomorphic
functional calculus and a 2-summable Fredholm module
$(F, \chh\otimes\cs)$
for ${\mathcal B}^\Gamma_0$ with 
Chern character given by the cyclic 2-cocycle
$\tau_{c,\Gamma}(A_0, A_1, A_2)$ which is equal to
$$\int_{X_\Gamma\times X\times X}
\Phi(z,x,y)\varpi(z,x,y)k_0(z,x,r) k_1(x,y,r) k_2(y,z,r)\,dz\,dx\,dy,$$
$r\in\Omega$, 
where the operators $A_0, A_1, A_2$  are in ${\mathcal B}^\Gamma_0$,
and whose Schwartz kernels are $k_0,k_1,k_2$ respectively.
Here $\Phi(z,x,y) = \int_\Delta \omega_\hyp$ is the oriented 
hyperbolic area of a geodesic triangle $\Delta$ with vertices at $x, y, z$.  
Furthermore
if $P(r)$ is a projection into a gap in the spectrum
of the Hamiltonian $H_{\eta, V}$. Then
    $P(r)$
    lies in a 2-summable dense subalgebra ${\mathcal B}^\Gamma_0$
of ${\mathcal B}^\Gamma$ and for almost any $r\in\Omega$ one has
$$   \Index(P(r)FP(r)) =   \langle\tau_{c,\Gamma}, [P(r)]\rangle
\in 2(g-1)\mathbb{Z}.$$
\end{thm}

\subsection{Summability of the Fredholm module}
The technical parts of the proof of the previous theorem
    rest on a lengthy calculation together with a key
estimate on kernels $k(z,w,r)$ on $\hyp \times \hyp\times \Omega$
which represent smooth functions of
the resolvent of $H+V$. This estimate has the form
$$ |k(z,w,r)|^2\leq C_2\exp(-C_3d(z,w)^2), \eqno(**)$$
where $C_2,C_3$ are constants (note that the RHS is independent of $r$).
This estimate is a result of \cite{BrSu}.
Since operators with kernels which have
    support in a band around the diagonal are
dense in the algebra  $\ck^\Gamma$
so too is the set of operators with kernels satisfying (**).
We denote by ${\mathcal B}^\Gamma_0$
the subalgebra
consisting of operators $A\in {\mathcal B}^\Gamma$,
with $[F,A]$ a Hilbert-Schmidt operator.
Now ${\mathcal B}^\Gamma_0$
is dense and
by \cite{Co} ${\mathcal B}^\Gamma_0$ is stable under
    the holomorphic functional calculus.
The last claim of the corollary on the range of values
taken by the cyclic cocycle follows using Morita equivalence
with $C(\Omega) \rtimes_{\bar\sigma} \Gamma$.
   The details are in \cite{Ma}\cite{CHMM}.

\subsection{The hyperbolic Connes-Kubo formula, part II}
We now have many cyclic 2-cocycles associated to our
model. We combine the cyclic 2-cocycles of subsection 1.4 
to produce a Connes-Kubo cocycle for the hyperbolic 
Hall conductance in Proposition \ref{prop:Connes-Kubo}, and our goal is to show 
that it is cohomologous to the Chern character of the
Fredholm module  $\tau_{c,\Gamma}$ as given in Theorem \ref{thm:Fred}.

%as follows.
%Given operators $A_0, A_1, A_2$ in $\ck^\Gamma$
%whose kernels $k_0,k_1,k_2$ are exponentially decaying
%    (cf equation (**) of the previous
%subsection) introduce the cyclic 2-cocycle $\tau_K$ defined by
%$$
%\tau_K(A_0, A_1, A_2) = \sum_{j=1}^g c_{j,j+g}(A_0, A_1, A_2)
%$$
%$$
%= \sum_{j=1}^{g}\int_{X_\Gamma\times X\times X}\varpi(z,x,y)
%\Psi_j(z,x,y)k_0(z,x,r)k_1(x,y,r)k_2(y,z,r)\,dz\,dx\,dy,
%$$
For $j=1, \ldots, g$, consider $\Psi_j(z,x,y)$ which is given by,
$$(\Omega_j(x)-\Omega_{j}(y))(\Omega_{j+g}(y)-\Omega_{j+g}(z))-
(\Omega_{j+g}(x)-\Omega_{j+g}(y))(\Omega_{j}(y)-\Omega_{j}(z)).$$
%We call this the Connes-Kubo cocycle and show now
%(following \cite{CHMM}) that it is cohomologous to $\tau_{c,\Gamma}$.
We claim first that
$\sum_{j=1}^g\Psi_j(z,x,y)$ is proportional to the `symplectic area'
of a triangle in $\real^{2g}$ with vertices $\Xi(x),\Xi(y),\Xi(z)$.
To prove this it suffices
to assume that the base point in $\hyp$ is one of the vertices of
the triangle, say $z$. Consider the expression
$$\sum_{j=1}^g\Psi_j(z,x,y)=\sum_{j=1}^g(\Omega_j(x)\Omega_{j+g}(y)-\Omega_{j+g}
(x)\Omega_{j}(y)).$$
Let $s$ denote the symplectic form on $\real^{2g}$ given by:
$s(u,v)=\sum_{j=1}^g(u_jv_{j+g} -u_{j+g}v_j).$
The `symplectic area' of a triangle $\Delta_E$
with vertices $0,\Xi(x),\Xi(y)$ is given by
$s(\Xi(x),\Xi(y))$. To appreciate why this is so
we need an argument from
{\cite{GH} (pp 333-336)}.
The form $s$ is the two form on $\real^{2g}$ given by
$$\omega_J=\sum_{j=1}^g du_j\wedge du_{j+g}.$$
Now the  symplectic area of a triangle $\Delta_E$ in $\real^{2g}$ with
vertices $0,\Xi(x),\Xi(y)$
is by definition the integral of $\omega_J$ over the triangle.
A brief calculation reveals that this yields
$s(\Xi(x),\Xi(y))/2$, proving our claim.
We have now established the following result.

\begin{prop}\label{prop:Connes-Kubo} 
The higher genus analogue of the Connes-Kubo formula is given by
the cyclic
2-cocycle $\tau_K$ on $\ck^\Gamma$ defined by
$$
\tau_K(A_0, A_1, A_2) = \sum_{j=1}^g \kappa\, c_{j,j+g}(A_0, A_1, A_2)
$$
$$
=\sum_{j=1}^{g}\int_{X_\Gamma\times X\times X} \kappa
\Psi_j(z,x,y)\varpi(z,x,y)k_0(z,x,r)k_1(x,y,r)k_2(y,z,r)\,dz\,dx\,dy
$$
for $r\in\Omega$.
Here the $k_j$ are the kernels of the $A_j, j=0,1,2$
(three exponentially decaying elements of $\ck^\Gamma$) and
$\sum_{j=1}^{g}\Psi_j(z,x,y)$ is proportional to the
`symplectic area' of the Euclidean triangle
$\Delta_E$ in $\real^{2g}$ with vertices
    $\Xi(x),\Xi(y),\Xi(z)$. Here $\kappa = 4\pi (g-1)/g$ is a constant
    depending only on the genus $g$, where $g>1$.   
    \end{prop}

%In the following argument we suppress the $r\in\Omega$ dependence
%as it is not important.
The constant $\kappa = 4\pi (g-1)/g$  is justified  in the discussion following Theorem~\ref{thm:Jost}
below. 
To compare the conductance cocycle $\tau_K$ with the Chern character 
cocycle $\tau_{c,\Gamma}$, we begin by recalling 
the following Theorem 5.5.1, page 222 in \cite{Jost}.${}^*$\footnote{${}^*$ Note that in \cite{CHMM} 
page 652, we used a different, incorrect argument
at this point, and we thank Siye Wu for pointing this out to us.}

\begin{thm}\label{thm:Jost}
Let $\Sigma$ be a compact Riemann surface of genus $g\ge2$ and $\alpha_1, 
\ldots, \alpha_g$ be a basis of holomorphic 1-forms on $\Sigma$. Then
$
\sum_{j=1}^g \alpha_j \otimes \bar\alpha_j
$
defines a K\"ahler metric on $\Sigma$ called the Bergman metric or the canonical metric,
that has nonpositive curvature vanishing at most at a finite number of points 
on $\Sigma$.
\end{thm}

It follows from this theorem, which uses the Riemann-Roch theorem,
 that $\omega_\alpha = 
\frac{\sqrt{-1}}{2}\sum_{j=1}^g \alpha_j \wedge \bar\alpha_j$
is a volume form on $\Sigma$. 
This is a subtle result as the holomorphic 1-form $\alpha_j$ cannot
be nowhere zero, which follows by an application of the 
Hopf index theorem, where we observe that the Euler characteristic is nonzero.
Therefore each term $\frac{\sqrt{-1}}{2} \alpha_j \wedge \bar\alpha_j$ 
by itself cannot be a volume form on $\Sigma$!

Next we recall the following basic fact relating holomorphic 1-forms
and harmonic 1-forms  on $\Sigma$. A 
(complex valued) 1-form $\alpha$ on $\Sigma$ is 
holomorphic if and only if $\alpha = a + \sqrt{-1} *a$, 
where $a$ is a (real valued) harmonic 1-form on $\Sigma$
and $*a$ is the Hodge $*$ of $a$.

If $a_j, \, j=1, \ldots 2g$ is a symplectic basis of harmonic 1-forms on $\Sigma$,
where $a_{j+g} = * a_j, \, j=1,  \ldots g$. Then $\alpha_j = a_j + \sqrt{-1} a_{j+g}$
is a basis of holomorphic 1-forms on $\Sigma$. By Theorem \ref{thm:Jost}
and its consequence, we deduce that $\sum_{j=1}^g a_j \wedge a_{j+g}$ 
is a volume form on $\Sigma$. 

Now let $\omega_\Sigma$ 
denote the volume form on $\Sigma = \hyp/\Gamma$ 
induced by the hyperbolic volume form $\omega_\hyp$
on $\hyp$. Then there is a positive constant $\kappa$ such that 
$\omega_\Sigma$ and $\kappa\sum_{j=1}^g a_j \wedge a_{j+g}$ are cohomologous.
%For simplicity, we rescale, $a_j \mapsto \frac{1}{\sqrt{\kappa}}a_{j}$,
%where we observe that since $\kappa$ is a constant, the rescaled 
%harmonic 1-forms are also harmonic. In this way, we can assume 
%without loss of generality that $\kappa =1$. 
To determine the constant $\kappa$, we integrate over the surface $\Sigma$
to get 
$$
\int_\Sigma \omega_\Sigma = \kappa \int_\Sigma \sum_{j=1}^g a_j \wedge a_{j+g}.
$$
Now each term $\int_\Sigma a_j \wedge a_{j+g} = 1$ by our choice of normalized 
symplectic basis.
By the Gauss-Bonnet theorem $\int_\Sigma \omega_\Sigma = 4\pi(g-1)$. 
Therefore $\kappa = 4\pi(g-1)/g$.

Thus by the argument above and Lemma \ref{lem:Jac}, we see that the difference
$\; \omega_\hyp - \kappa \Xi^*(\omega_J) = d\Lambda, \;$ where $\Lambda$
is a $\Gamma$-invariant 1-form on $\hyp$.
More particularly
for a  geodesic triangle $\Delta\subset\hyp$ with vertices at $x,y,z \in \hyp$,
$$
\begin{array}{lcl}
\displaystyle\int_\Delta\omega_\hyp
&=&\displaystyle \kappa \int_\Delta\Xi^*(\omega_J) + 
\int_\Delta d\Lambda\\[+7pt]
&=&\displaystyle \kappa \int_{\Xi(\Delta)}\omega_J + \int_{\partial\Delta} \Lambda
\end{array}
$$
Now $\Xi$ cannot map geodesic triangles
to Euclidean triangles in $\real^{2g}$
as $\Xi({\Delta})$ is a compact subset of a non-flat embedded two
dimensional surface in $\real^{2g}$.
Moreover as  $\Psi_j(z,x,y)=0$ whenever the images of $z,x,y$ under $\Xi$
lie in a Lagrangian subspace (with respect to the symplectic
form $s$) of $\real^{2g}$, $\tau_K$ and $\tau_{c,\Gamma}$
are not obviously proportional.

%To see that (suitably normalised) they are cohomologous we
%first renormalise $\omega_J$ so that
%$\Xi^*(\omega_J)=\omega_\hyp$
Next we write $\omega_J= d\theta$.
Considering the difference $\tau_K-\tau_{c,\Gamma}$ one sees that the
key is to understand
$$\int_{\Xi(\Delta)}\omega_J\quad -\int_{\Delta_E}\omega_J
=\int_{\partial\Xi(\Delta)}
\theta \quad - \int_{\partial\Delta_E}\theta.
$$
Now this difference of integrals around the boundary can be written as
the sum of three terms corresponding to splitting the boundaries
$\partial\Xi(\Delta)$ and $\partial\Delta_E$ into
three arc segments each. We introduce some notation for this,
writing
$$\partial\Xi(\Delta)=\Xi(\ell(x,y))\cup\Xi(\ell(y,z))\cup\Xi(\ell(z,x)),$$
where $\ell(x,y)$ is the geodesic in $\hyp$ joining $x$ and $y$
(with the obvious similar definition of the other terms). We also write
$$\partial\Delta_E=m(x,y)\cup m(y,z)\cup m(z,x),$$
where $m(x,y)$ is the straight line joining $\Xi(x)$ and $\Xi(y)$
(and again the obvious definition of the other terms).
Then we have
$$\int_{\partial\Xi(\Delta)}\theta\quad - \int_{\partial\Delta_E}\theta=
h(x,y)+h(y,z)+h(z,x) \eqno(*)$$
where $h(x,y)= \int_{\Xi(\ell(x,y))}\theta - \int_{m(x,y)}\theta$
with similar definitions  for $h(y,z)$ and $h(z,x)$.

Notice that we have $h(x,y)=\int_{D_{xy}}\omega_J$
where $D_{xy}$ is a disc with boundary $m(x,y)\cup \Xi(\ell(x,y))$. From
    this it is easy to see that $h(\gamma x,\gamma y)=h(x,y)$
for $\gamma\in \Gamma$.

Now consider $j(x,y) = \int_{\ell(x,y)} \Lambda$. 
Since $\Lambda$ is $\Gamma$-invariant, it follows that
$j(\gamma x,\gamma y)=j(x,y)$
for $\gamma\in \Gamma$.
Then by the 
computation done above, we see that
$$
\displaystyle\int_\Delta\omega_\hyp = 
\kappa \int_{\Delta_E}\omega_J + \kappa (h(x,y)+h(y,z)+h(z,x))
+ j(x,y)+j(y,z)+j(z,x)
$$
%and then 
We normalise $\sum_{j=1}^g\Psi_j(z,x,y)$ so that it equals
$\int_{\Delta_E}\omega_J$.
Then, 
$$
 \Phi(x,y,z) = \kappa \sum_{j=1}^g\Psi_j(z,x,y) + \partial(\kappa  h+j)(x,y,z)
$$
where $ \Phi(x,y,z)  = \displaystyle \int_\Delta\omega_\hyp. $

Introduce the bilinear functional $\tau_1$ on $\ck^\Gamma$
given by
$$
\begin{array}{lcl}
\tau_1(A_0,A_1) & = & \int_{X_\Gamma \times X}
(h(x,y) + j(x,y))k_0(x,y)k_1(y,x)\,dx\,dy \\[+7pt]
& = &\tr_{\ck^\Gamma}(A_{\kappa h + j}A_1),
\end{array}
$$
where the operator $A_j$ has kernel $k_j(x,y, r),\;j=0,1$ 
and $A_{\kappa h + j}$ 
is the operator
with kernel $(\kappa h(x,y) + j(x,y)) k_0(x,y, r)$.
So we have proved that formally the two cyclic 2-cocycles 
satisfy,
$$
b\tau_1= \tau_K-\tau_{c,\Gamma},
$$
where $b$ is the Hochschild boundary operator, so that 
they are cohomologous cyclic 2-cocycles. What remains is 
to understand the domain of the cochains, which is what is 
addressed next.

We want to see that $\tau_1$ is densely defined. By Theorem 1.5, one has
an isomorphism
$$
\Phi_F: {\mathcal B}^\Gamma\cong C(\Omega) \rtimes_{\bar\sigma} 
\Gamma\otimes {\mathcal K}
(L^2(F)).
$$
Here $F$ denotes a fundamental domain for the action of $\Gamma_g$
on $\mathbb H$.
Now any element $x$ in $C(\Omega) \rtimes_{\bar\sigma} \Gamma\otimes 
{\mathcal K}$
can be written as a matrix $(x_{ij})$, where $x_{ij} \in C(\Omega)
\rtimes_{\bar\sigma} \Gamma$.
So we can define
$$N_k(x) = ( \sum_{i,j}\nu(x_{ij})^2)^{\frac{1}{2}},$$
where $$\nu(x_{ij}) = (\sum_{h\in \Gamma_g} (1+\ell(h)^{2k})
|x(h)|^2)^{\frac{1}{2}}$$
and $\ell$ denotes the word length function on the group $\Gamma_g$.
Using a slight modification of the argument given in \cite{Co2},
III.5.$\gamma$, one can prove that there is a subalgebra
${\mathcal B}^\Gamma_\infty$ of ${\mathcal B}^\Gamma$ which\\
(i) contains
$C(\Omega) \rtimes_{\bar\sigma, alg} \Gamma\otimes {\mathcal R}$,
where $\mathcal R$
denotes the algebra of smoothing operators on $F$ and 
$\rtimes_{\bar\sigma, alg}$
denotes the algebraic twisted crossed product,\\
(ii) is stable
under the holomorphic functional calculus, and \\
(iii) is such that $N_k(x)<\infty$
for all $x\in {\mathcal B}^\Gamma_\infty$ and $k\in \mathbb N$.

Then, following \cite{Co2}, we have that
   the trace $\tau \otimes \mbox{Tr}$ on
$C(\Omega) \rtimes_{\bar\sigma, alg} \Gamma\otimes {\mathcal R}$, is continuous
for the norm $N_k$, for $k$ sufficiently large, and thus extends by continuity
to ${\mathcal B}^\Gamma_\infty$. Note that elements in
${\mathcal B}^\Gamma_\infty$ have Schwartz kernels which have
rapid decay away from the diagonal. The next result summarises the
discussion above.

\begin{prop} The algebra ${\mathcal B}^\Gamma_\infty$ is dense
in ${\mathcal B}^\Gamma$, is closed under the
holomorphic functional calculus
and is contained in the ideal ${\mathcal I}$ of
$\ck^\Gamma$ consisting of operators with finite trace.
\end{prop}

Now $\tau_K$ is defined on  ${\mathcal B}^\Gamma_\infty$ while
$\tau_{c,\Gamma}$ is defined on  ${\mathcal B}^\Gamma_0$
as we noted earlier. Both of these algebras contain the
operators whose Schwartz kernels are supported in a band around the
diagonal. Thus the subalgebra
${\mathcal B}^\Gamma_\infty\cap {\mathcal B}^\Gamma_0$ is dense
and stable under the holomorphic functional calculus.
%If $b$ denotes the Hochschild coboundary map then
%(*) says that $b\tau_1= \tau_K-\tau_{c,\Gamma}$.
%The Lipschitz property of the Jacobi map means that 
Since $\Lambda$ is $\Gamma$-invariant, it is bounded, therefore
$|j(x,y)| \le ||\Lambda||_\infty d(x,y)$, where $||\Lambda||_\infty$ is the supremum
norm of $\Lambda$ and $d(x,y)$ is  the hyperbolic distance from $x$ to $y$.
An explicit expression for $\theta$ shows that it grows linearly in terms of
$d(x,y)$, so that $h(x,y)$  grows at worst like $d(x,y)^2$. (for more details, see \cite{CHMM})
Therefore if $A_0\in\ck^\Gamma_\infty$ then so too does $A_{\kappa h +j}$.
Hence we have $\tau_1$ defined on $\ck^\Gamma_\infty\cap{\mathcal B}^\Gamma_0$.
This section has proved our main theorem.

\begin{thm}
The Connes-Kubo cocycle $\tau_K$ and the Chern character
cocycle $\tau_{c,\Gamma}$ arising
as the Chern class of the
Fredholm module $(F, \chh\otimes\cs)$, are cohomologous
as cyclic cocycles on $\ck^\Gamma_\infty\cap{\mathcal B}^\Gamma_0$.
\end{thm}

%%%%%%%%%%%%%%%%%%%%%%%%%%%%%
\section{Appendix : On the Quantum Adiabatic Theorem (QAT)}
%%%%%%%%%%%%%%%%%%%%%%%%%%%%%%

One knows that the (time) evolution determined by a time
{\em independent} Hamiltonian reduces to the spectral
theory of the Hamiltonian. The
QAT says that the (time) evolution of a slowly
varying time {\em dependent} Hamiltonian reduces to the spectral
theory of an associated family of adiabatic Hamiltonians.
The setting for the QAT is as follows.
Let $s\to H(s)$ be a smooth family of Hamiltonians (self-adjoint
operators) $\tau=\text{time scale}$ and $s=t/\tau=\text{scaled
time}$. Consider now
the {\em physical evolution}
\[
i\partial_t U(t)=H(t/\tau)U(t),\quad U(0)=1
\]
or equivalently
\begin{equation}
\label{eq1}
i\partial_s U_{\tau} (s)=\tau H_{\tau}(s) U_{\tau}(s),\quad U_{\tau}(0)=1\,.
\end{equation}
Let $P(0)$ denote the spectral projection onto a gap in the spectrum of $H(0)$,
that is we have
$P(0)=\chi_{(-\infty,E]}(H(0))$
where $E\not\in$ spectrum of $(H(0))$.

The {\em adiabatic evolution} is determined by the equation
\begin{equation}
\label{eq2}
P(s)=U_a(s) P(0) U_a(s)^*,\quad U_a(0)=1
\end{equation}
where $P(s)$ denotes spectral projection onto a gap in the spectrum
of $H(s)$. Let $H_a(s)$ denote the generator of $U_a(s)$. It is
also known as the {\em adiabatic Hamiltonian} and is given
by
\begin{equation}
\label{eq3}
H_a(s)=\frac{i}{\tau} (\partial_s U_a(s)) U_a(s)^*
\end{equation}
%\end{setting}
\begin{lem} The adiabatic Hamiltonian
$H_a(s)$ satisfies the {equation of motion}
\[
[H_a(s), P(s)]=\frac{i}{\tau} \partial_s P(s)
\]
\end{lem}
{\em Proof.}
Differentiating (\ref{eq2}), we have
\begin{align*}
\partial_s P(s) &= \partial_s U_a(s) P(0) U_a(s)^*+U_a(s) P(0) 
\partial_s U_a(s)^*\\
     &=\partial_s U_a(s) P(0) U_a(s)^*-U_a(s) P(0) U_a(s)^* \partial_s U_a(s)
       U_a(s)^* \\
     &=(\partial_s U_a(s))U_a(s)^* U_a(s) P(0) U_a(s)^*-\frac{\tau}{i}
P(s) H_a(s)\\
     &=\frac{\tau}{i} [H_a(s), P(s)]
\end{align*}
${}$\hfill
%\end{proof}
\begin{lem} Let $f$ be a measurable function on $\mathbb R$. Then
$\ds H_a(s)=f(H(s))+\frac{i}{\tau}[\partial_sP(s),P(s)]$
satisfies the equations of motion.
\end{lem}
{\em Proof.} $[f(H(s)),P(s)]\equiv 0$ and
$\left[[\partial_s P(s), P(s)], P(s) \right] =\partial_s P(s)$
since $P(s)^2=P(s)$ and $P(s)$ is a spectral projection of $H(s)$.
Define the adiabatic Hamiltonian as
\begin{equation}
\label{eq4}
H_a (s) =H(s)+\frac{i}{\tau} [\partial_s P(s), P(s)]
\end{equation}
Then equation (\ref{eq2}) is satisfied and
$U_a(s):\; \operatorname{Range}P(0)\rightarrow\operatorname{Range}P(s)$
i.e. the initial value problem
\[
i\partial_s \psi(s)=\tau H_a(s) \psi(s), \quad \psi(0)\in
\operatorname{Range}P(0)
\]
has the property that $\psi(s)\in\operatorname{Range} P(s)$ $\forall s$.
%\end{proof}
\begin{thm}[Quantum Adiabatic Theorem (QAT) \cite{Av+S+Y}]
Let $s\to H(s)$  be a smooth family of self-adjoint Hamiltonians
and $s\to P(s)$ be a smooth family of spectral projections as before
such that
\[
\sup \{\| P(s)\|<\infty \mid s\in [0,\infty)\}
\]
and the commutator equation
$[\partial_s P(s), P(s)]=[H(s),X(s)]$
has an operator-valued solution
$X(s)$, such that \!$X(s)$ and $\partial_s X(s)$ are {\em bounded}. Then
one has
\[
   \|(U_{\tau} (s) - U_a (s))P(0)\|\le \frac{1}{\tau}
       \max\limits_{s\in [0,\infty)}\{2\|X(s) P(s)\|+
       \|\partial_s(X(s) P(s)) P(s)\|\}
\]
That is, the adiabatic evolution $U_a(s)$ approximates the physical
evolution $U_{\tau}(s)$ as the adiabatic parameter $\tau\to \infty$.
Equivalently, the adiabatic Hamiltonian $H_a(s)$ approximates the
physical  Hamiltonian $H_{\tau}(s)$ on the range of $P$,
as the adiabatic parameter $\tau\to \infty$.
\end{thm}

Note that the 
hypotheses on $P(s)$ are satisfied if $P(s)$ is a spectral projection
onto a gap in the spectrum of $H(s)$ because one can then
define
\[
X(s)=\frac{1}{2\pi i} \oint_C R(z,s) \partial_s P(s) R(z,s) dz
\]
where $C$ is a contour in $\mathbb{C}$ enclosing the spectrum
in $(-\infty,E]$, $E\not\in$ $\operatorname{spec}$ $(H(s))$ and
$R(z,s)=(H(s)-z)^{-1}$ is the resolvent.

%%%%%%%%%%%%%%%%%%%%%%%%%%%%
\section {Appendix: Conductance cocycles}
%%%%%%%%%%%%%%%%%%%%%%%%%%%

In this subsection we present  an argument which
derives from physical principles the hyperbolic Connes-Kubo formula for the
`Hall conductance'.
Our reasoning is that
   the Hall conductance in the Euclidean situation
is measured experimentally by determining the equilibrium
   ratio of the current in the direction of the applied electric field to the
   Hall voltage, which is the potential difference in the orthogonal direction.
   To calculate this mathematically we instead determine
   the component of the induced current that is orthogonal to the applied
   potential. The conductance can then be obtained by dividing this
   quantity by the magnitude of the applied field.
In the hyperbolic case  preferred directions.
are obtained by interpreting the generators of the fundamental group
as geodesics on hyperbolic space giving a family of preferred directions
emanating from the base point. For each pair of
directions it is therefore natural to imitate the procedure of the
Euclidean case
and mathematically this is done as follows.

The Hamiltonian $H$ in a magnetic field depends on the magnetic vector
potential ${\bf A}$ and the functional derivative $\delta_kH$
of $H$ with respect to one of the components of ${\mathbf A}$, denoted
$A_k$, gives
   the current density $J_k$,
where we consider adiabatic variations within a one-parameter family
$A_k(s)$, which we can choose without loss of generality to be bounded,
since ${\bf A}(0) = - \theta \frac{dx}{y}$ defines
a  bounded operator in the hyperbolic metric.
The expected value of the current in a state described by a projection
operator $P$ into a spectral gap of $H$ is therefore $\tr(P\delta_kH)$
(cf \cite{Av+S+Y} equation (3.2)). (Note that an argument,
using the fact that $P$ is a member of a family $P(s)$
of projections which correspond to gaps for small $s$,
is required to see that $P\delta_kH$ is trace class.)
The following lemma is not proved by a rigorous argument: one needs to
check various analytical details as in \cite{Xia} which we omit as they
would take us too far afield.
For this discussion $\tr$ will denote a generic trace.
\begin{lem} In the adiabatic limit as the adiabatic parameter
$\tau\to \infty$, the functional derivative of the adiabatic Hamiltonian
$\delta_k H_a(s)$ approximates the functional derivative of the
physical  Hamiltonian $\delta_k H_{\tau}(s)$ on the range of $P$, and
one has
$$\tr(P\delta_kH) = i\tr(P[\partial_t P,\delta_kP]).$$
\end{lem}
{\em Proof.}
The first statement is the result of a calculation.
It uses the explicit forms of $\delta_k$ and $H_a(s)$
and the fact that the family $A_k(s)$ is bounded
to show that the norm of the difference $\delta_k H_a(s)-\delta_k H_{\tau}(s)$
goes to zero as $s\to 0$.
By using the invariance of the trace under the adjoint action of operators
and the equation of motion we see that
\begin{align*}
\tr(P[\partial_tP,\delta_kP]) &= -\tr([P,\delta_kP]\partial_tP)\\
&= -i\tr([P,\delta_kP][P,H_a])\\
&= i\tr([P,[P,\delta_kP]]H_a).
\end{align*}
Now $\delta_kP = \delta_k(P^2) = P(\delta_kP) + (\delta_kP)P$, whence
$P(\delta_kP) P = 0$ and we have
\begin{align*}
[P,[P,\delta_kP]] &= P(P(\delta_kP) - (\delta_kP)P)
-(P(\delta_kP) - (\delta_kP)P)P\\
&= P(\delta_kP) +(\delta_kP)P = \delta_kP.
\end{align*}
Consequently we may write
$$\tr(P[\partial_tP,\delta_kP]) = i\tr((\delta_kP)H_a)
= i\tr(\delta_k(PH_a)) - i\tr(P(\delta_kH_a)),$$
and, assuming that the trace is invariant under variation of $A_k$, the
first term vanishes.  The result asserted follows by taking the limit
as the adiabatic parameter
$\tau\to \infty$.
%\end{proof}
By following \cite{KS}, one sees that in fact the limit of the lemma
is true to all orders. We note further
that if the only $t$-dependence in $H$ and $P$ is due to the
adiabatic variation of $A_j$, a
component distinct from $A_k$, then
$\partial_t = \partial A_j/\partial t\times\delta_j$.
Working  in the Landau gauge so that the electrostatic
potential vanishes, the electric field is given by
${\bf E} = -\partial {\bf A}/\partial t$, and so $\partial_t = -E_j\delta_j$.
Combining this with the previous argument we arrive at the following result:
\begin{cor}
The conductance for currents in the $k$ direction induced by electric fields
in the $j$ direction is given by
$-i\tr(P[\delta_jP,\delta_kP])$.
\end{cor}
{\em Proof.}
The expectation of the current $J_k$ is given by
$$\tr(P\delta_kH) = i\tr(P[\partial_tP,\delta_kP])
= -iE_j\tr(P[\delta_jP,\delta_kP]),$$
from which the result follows immediately.

%%%%%%%%%%%%%%%%%%%%%%%%%%%%%

%%%%%%%%%%%%%%%%%%%%%%%%%%%%%%%%

\end{document}